\documentclass[12pt]{amsart}

\usepackage{amssymb,latexsym,amsmath,extarrows, mathrsfs, amsthm}
\usepackage{graphicx}

\usepackage[margin=1in, centering]{geometry}

\newtheorem{theorem}{Theorem}[section]
\newtheorem{lemma}[theorem]{Lemma}

\newtheorem{definition}[theorem]{Definition}
\newtheorem{corollary}[theorem]{Corollary}

\newcommand{\cs}{\mathcal S}

\newcommand{\nf}{\infty}

\newcommand{\ZR}{\mathbb{R}}
\newcommand{\ZZ}{\mathbb{Z}}
\newcommand{\ZN}{\mathbb{N}}

\newcommand{\bT}{{\bf T}}

\newcommand{\hichi}{\raisebox{0.7ex}{\(\chi\)}}

\begin{document}

\title[transforms along arithmetic functions]{Discrete bilinear Radon transforms along arithmetic functions with many common values}
\author{Dong Dong}
\address{Mathematics Department\\
University of Illinois at Urbana-Champaign\\
1409 West Green Street, Urbana, IL 61801, USA}
\email{ddong3@illinois.edu}
\author{Xianchang Meng}
\address{Centre de Recherches Math\'{e}matiques\\
Universit\'{e} de Montr\'{e}al\\
P.O. Box 6128,
Centre-ville Station\\
Montr\'{e}al (Qu\'{e}bec)
H3C 3J7, Canada}
\email{meng@crm.umontreal.ca}
\date{\today}

\begin{abstract}
We prove that for a large class of functions $P$ and $Q$, there exists $d\in (0,1)$ such that the discrete bilinear Radon transform $$B^{\rm dis}_{P,Q}(f,g)(n)=\sum_{m\in\mathbb{Z}\setminus\{0\}} f(n-P(m))g(n-Q(m))\frac{1}{m}$$ is bounded from $l^2\times l^2$ into $l^{1+\epsilon}$ for any $\epsilon\in (d,1)$. In particular, the boundedness holds for any $\epsilon\in (0,1)$ when $P$ (or $Q$) is the Euler totient function $\phi(|m|)$ or the prime counting function $\pi(|m|)$. 
\end{abstract}
{\let\thefootnote\relax
\let\thefootnote\relax}

\maketitle
\section{Introduction}
\setcounter{equation}0
One of the most important multilinear operators in harmonic analysis is the bilinear Hilbert transform (BHT), which is defined by
\[
B(f,g)(x)=\int f(x-t)g(x+t)\, \frac{dt}{t}, \,\, f,g\in\cs(\ZR),
\]
where $\cs(\ZR^n)$, $n\in\ZN$, is the Schwartz space on $\ZR^n$. Lacey-Thiele's breakthrough study of BHT \cite{LT97,LT99} has generated numerous investigations in multilinear operators. One direction of generalizing BHT is to replace the linear terms $+t$ and $-t$ with some non-linear polynomials $P$ and $Q$, and consider the operator
\[
B_{P,Q}(f,g)(x):=\int f(x-P(t))g(x-Q(t))\, \frac{dt}{t}, \,\, f,g\in\cs(\ZR).
\]
See the articles of Li \cite{LNY,LAP}, Li-Xiao \cite{LX} and the first author \cite{D} for some recent results about $B_{P,Q}$. Motivated by many works on discrete linear operators (see, for example, Bourgain \cite{B}, Ionescu-Wainger \cite{IW}, Krause \cite{K}, Mirek-Trojan \cite{MT}, Mirek-Stein-Trojan \cite{MST2, MST1}, Pierce \cite{Lillian}, Zorin-Kranich \cite{ZK}, etc.), one may also consider a discrete analogue of BHT:
\[
B^{\text{dis}}(f,g)(n):=\sum_{m\in\ZZ\setminus\{0\}} f(n-m)g(n+m)\frac{1}{m}, \,\, f,g\in D(Z),
\]
or more generally,
\begin{equation} \label{def of B}
B_{P,Q}^{\text{dis}}(f,g)(n):=\sum_{m\in\ZZ\setminus\{0\}} f(n-P(m))g(n-Q(m))\frac{1}{m}, \,\, f,g\in D(Z).
\end{equation}
Here $D(\ZZ^n)$ is the space of compactly supported complex-valued functions defined on $\ZZ^n$. It turns out that the boundedness of $B^{\text{dis}}$ is equivalent to that of BHT by transference principle \cite{BMMS}. Therefore, Lacey-Thiele's Theorem already covers the boundedness of $B^{\text{dis}}$. However, transference is not available if $-m$ and $+m$ are replaced with non-linear functions. In fact, even the $l^2\times l^2\to l^{1,\nf}$-boundedness of $B_{P,Q}^{\text{dis}}$ for $P(m)=m$ and $Q(m)=m^2$ is an extremely difficult problem and still out of reach of current popular techniques such as time-frequency analysis and circle methods.

Hu and Li \cite{HL} first obtained the $l^2\times l^2\to l^{1+\epsilon}$ boundedness of $B_{P,Q}^{\text{dis}}$ when $P(m)=m$ and $Q(m)=m^2$. Recently, using a different method, the first author \cite{DCRM} extended Hu-Li's result to function pairs satisfying the \textbf{Condition $(\star)$}. More precisely, the following definition and theorem were given in \cite{DCRM}.
\begin{definition}
For any function $R$ that maps $\ZZ$ into $\ZZ$, let
\begin{equation} \label{def: APQ}
A^R:=\left\{(m_1,m_2)\in(\ZZ\setminus\{0\})^2: R(m_1)=R(m_2)\right\}.
\end{equation}
We say that $R$ satisfies \textbf{Condition $(\star)$} if there are constants $D_1$ and $D_2$ such that $\frac{|m_1|}{|m_2|}\le D_1$ for all $(m_1,m_2)\in A^R$, and for each $m_1\in\ZZ$ there are at most $D_2$ pairs $(m_1,m_2)$ in the set $A^R$.
\end{definition}
\begin{theorem}[\cite{DCRM}, Theorem 1.1  \footnote{In the original statement of this theorem, Condition $(\star)$ is imposed on $P-Q$ instead of $P$ (or $Q$), but the same arguments there in fact work in both cases. We are indebted to Xiumin Du for pointing this out.}]   \label{old thm}
Given two functions $P$ and $Q$ that map $\ZZ$ into $\ZZ$,
assume that $P$ or $Q$  satisfies Condition $(\star)$. Then for any $\epsilon\in(0,1]$, there is a constant $C_{\epsilon}$ depending only on $\epsilon$, $D_1$ and $D_2$ such that the operator $B_{P,Q}^{\text{dis}}$ defined by \eqref{def of B} satisfies
\begin{equation} \label{eq: main}
\|B_{P,Q}^{\text{dis}}(f,g)\|_{l^{1+\epsilon,\nf}}\le C_{\epsilon}\|f\|_{l^2}\|g\|_{l^2}, \text{ for any } f,g\in l^2.
\end{equation}
\end{theorem}

Monotonic functions and non-constant polynomials satisfy Condition $(\star)$. However, Condition $(\star)$ requires that the function can attain each value for only bounded number of times, which excludes numerous arithmetic functions. For example, Ford \cite{Ford} proves that for any $k\ge 2$ there exists $n_k$ such that the Euler's totient function $\phi(n)$ equals $n_k$ for at least $k$ times. Thus $\phi$ does not satisfy Condition $(\star)$. By the fact that gaps between primes can be arbitrarily large, the prime counting function $\pi(n)$ does not satisfy Condition ($\star$) either. Due to the discrete nature of the operator $B_{P,Q}^{\text{dis}}$, it is interesting to seek for a weaker condition for $P$ and $Q$ (under which \eqref{eq: main} still holds) that includes some important arithmetic functions having many common values. The definition and the main theorem of our paper below serve as this purpose.

\begin{definition}
For any function $R$ that maps $\ZZ$ into $\ZZ$, let
\[
S_{M,N}^R:=\left\{(m,n): R(m)=R(n), \frac{1}{2}N\le |n|\le 2N, \frac{1}{2}M\le |m|\le 2M\right\}. 
\]
We say that $R$ satisfies \textbf{Condition $(\star\star)$} if there exist constants $\delta>0$ and $\delta'>0$ such that
\begin{equation} \label{SMN}
|S_{M,N}^R|\leq \frac{\delta' MN}{(\log M\log N)^{1+\delta}}. 
\end{equation}
\end{definition}
Roughly speaking, Condition $(\star\star)$ says that the solutions of $R(m)=R(n)$ in each dyadic strip have density slightly less than that of prime pairs. It is easy to see that Condition $(\star)$ implies Condition $(\star\star)$ for any $\delta>0$ and thus the following main theorem of this paper extends Theorem \ref{old thm}.

\begin{theorem} \label{new thm}
Given two functions $P$ and $Q$ that map $\ZZ$ into $\ZZ$,
assume that $P$ or $Q$ satisfies Condition $(\star\star)$. Then for any $\epsilon\in(\frac{1}{2\delta+1},1)$, there exists a constant $C_{\delta,\delta'}$ depending only on $\delta$ and $\delta'$ appeared in Condition $(\star\star)$ such that the operator $B_{P,Q}^{\text{dis}}$ defined by \eqref{def of B} satisfies
\begin{equation}
\|B_{P,Q}^{\text{dis}}(f,g)\|_{l^{1+\epsilon}}\le C_{\delta, \delta'}\|f\|_{l^2}\|g\|_{l^2}, \text{ for any~} f,g\in l^2.
\end{equation}
\end{theorem}

Remarks. (1). When $\epsilon=\frac{1}{2\delta+1}$ or $\epsilon=1$, we have weak-$l^{1+\epsilon}$ estimate. See the proof of Theorem \ref{new thm} in Section \ref{section: proof}.

(2)  When $P(m)=m$ and $Q$ is a polynomial, the operator norm of $B_{P,Q}^{\text{dis}}$ we obtained is independent of $Q$. Such uniform estimates also appear in the continuous setting (\cite{GL,LRev,LX,Thiele}).

(3) Note that the lower bound for $\epsilon$ goes to $0$ as $\delta$ tends to $\nf$. This gives an evidence that $l^2\times l^2\to l^1$-boundedness of $B_{P,Q}^{\text{dis}}$ may be true for at least some special $P$ and $Q$. 


Very interestingly, Condition $(\star\star)$ covers some important arithmetic functions from number theory. 
\begin{corollary}\label{coro}
If $P$ or $Q$ equals the Euler's totient function $\phi(|m|)$ or the prime counting function $\pi(|m|)$, then for any $\epsilon\in (0, 1)$, we have
\begin{equation}
\|B_{P,Q}^{\text{dis}}(f,g)\|_{l^{1+\epsilon}}\le C_{\epsilon}\|f\|_{l^2}\|g\|_{l^2}, \text{ for any~} f,g\in l^2.
\end{equation}
\end{corollary}

\begin{figure}[h]
\centering
\begin{minipage}[c]{0.5\textwidth}
\centering
\includegraphics[width=0.9\textwidth]{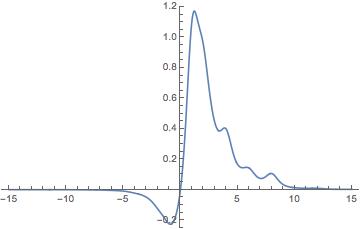}
\caption{\footnotesize $B_{P, Q}^{\rm dis}(f,g)(x)$: $P(m)=\phi(|m|)$, $Q(m)={\rm sgn}(m)d(|m|)$}
\label{fig-Euler-n}
\end{minipage}%
\begin{minipage}[c]{0.5\textwidth}
\centering
\includegraphics[width=0.9\textwidth]{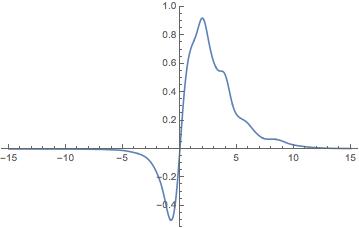}
\caption{\footnotesize $B_{P, Q}^{\rm dis}(f,g)(x)$: $P(m)=\pi(|m|)$, $Q(m)={\rm sgn}(m)d(|m|)$}
\label{fig-Pi-n}
\end{minipage}%
\end{figure}

\bigskip
To better demonstrate the behaviors of $B_{P,Q}^{\text{dis}}(f,g)(n)$ when $P(m)=\phi(|m|)$ or $P(m)=\pi(|m|)$, we exhibit the graphs of the operator for $f(x)=g(x)=\frac{1}{x^2+1}$. Fix $Q(m)={\rm sgn}(m)d(|m|)$, where $d(m):=\sum_{a|m} 1 ~(m>0)$ is the divisor function. We have 
$$
|S_{M,N}^{d}|\gtrsim \frac{MN}{\log M\log N},
$$
as $d(p)=2$ for any prime $p$. Therefore $Q$ does not satisfies Condition $(\star\star)$. We truncate the sum $B_{P, Q}^{\rm dis}(f,g)(x)=\sum_{m\in\ZZ\setminus\{0\}} f(x-P(m))g(x-Q(m))\frac{1}{m}$ up to $|m|\leq T_0=1000$. For any $|x|\leq 15$, the error will be bounded by $\frac{1}{(T_0-15)^2}$ which is good enough for us to plot Figure \ref{fig-Euler-n} and \ref{fig-Pi-n}.

\noindent

\subsection{Open Problems}

There are a few related open problems to consider.

(1). When $P$ and $Q$ are polynomials, we believe that the operator norm of $B_{P,Q}$ may be chosen to be independent of the coefficients of both $P$ and $Q$. We do not know how to achieve this. 

(2). A useful operator related with $B_{P,Q}^{\text{dis}}$ is the corresponding maximal operator $$B_{P,Q}^*(f,g)(n)=\sup_{M\in[1,\nf)}\left|\frac{1}{M}\sum_{m=1}^M f(n-P(m))g(n-Q(m))\right|.$$ It is conjectured that this operator is bounded from $l^2\times l^2$ into $l^{1,\nf}$. See \cite{DCRM,HL} for some positive results about this operator.

(3). In our proof of Corollary \ref{coro}, we use results about the Carmichael conjecture and gaps between primes. In converse, we wonder if the boundedness of this kind of operators could imply some information of the Carmichael conjecture and prime gaps. 

(4). Note that if $P$ is a constant function, then $B_{P,Q}^{\text{dis}}$ is bounded using the theory of discrete linear Radon transform \cite{IW,M}. In this extreme case, $|S_{M,N}^{P}|\simeq MN$. It remains to understand what happens if $|S_{M,N}^{P}|$ lies in between: 
\begin{equation} \label{not satisfy condition double star}
\frac{MN}{\log M\log N}\lesssim |S_{M,N}^{P}|\lesssim MN.
\end{equation}
For example, besides the divisor function $d$,   M\"{o}bius function $\mu$ and $\Omega$ function (the number of prime divisors) also satisfy \eqref{not satisfy condition double star}. Therefore, Theorem \ref{new thm} does not cover the cases when both $P$ and $Q$ are among these functions. Nevertheless, let us examine the graphs of the operator as before: 
\noindent\begin{figure}[h]
\centering
\begin{minipage}[c]{0.5\textwidth}
\centering
\includegraphics[width=0.9\textwidth]{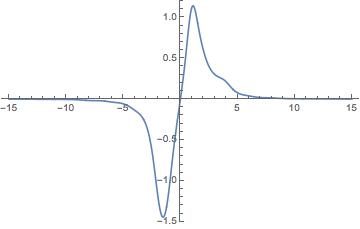}
\caption{\footnotesize $B_{P, Q}^{\rm dis}(f,g)(x)$: $P(m)=\mu(|m|)$, $Q(m)={\rm sgn}(m)d(|m|)$}
\label{fig-Divisor-n}
\end{minipage}%
\begin{minipage}[c]{0.5\textwidth}
\centering
\includegraphics[width=0.9\textwidth]{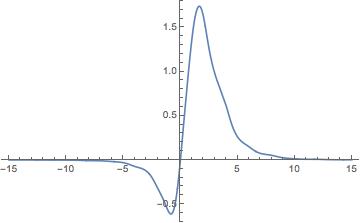}
\caption{\footnotesize $B_{P, Q}^{\rm dis}(f,g)(x)$: $P(m)=\Omega(|m|)$, $Q(m)={\rm sgn}(m)d(|m|)$}
\label{fig-Divisor-Omega}
\end{minipage}%
\end{figure}

These pictures have similar shapes as those in Figures \ref{fig-Euler-n}-\ref{fig-Pi-n}. It is reasonable to conjecture that our main theorem still holds in these cases.

\bigskip 


Throughout this paper, we use $A\lesssim B$ to denote the statement that $A\le CB$ for some positive constant C. When the implied constant $C$ depends on some parameter, say $\delta$, we may write $A\lesssim_{\delta} B$. $A\simeq B$ is short for $A\lesssim B$ and $B\lesssim A$. For any set of integers $E$, $|E|$ and $\hichi_E$ will be used to denote the counting measure and the indicator function of $E$, receptively.

\section{Arithmetic Functions with Many Common Values} \label{section: arithmetic functions}
\setcounter{equation}0
In this section, we show that Euler's totient function $\phi(n)$ and the prime counting function $\pi(n)$ satisfy Condition $(\star\star)$, and thus prove Corollary \ref{coro} assuming Theorem \ref{new thm}. We will introduce some backgrounds in a friendly way, as to make our paper more readable to both number theorists and analysts. 

\subsection{Euler's totient function}
Euler introduced the function $\phi(n)$ which counts the number of positive integers $\leq n$ that are coprime to $n$. Euler's totient function not only has deep connections with prime numbers, but also appears in many classical theorems in number theory. We know that Euler's totient function $\phi(n)$ is multiplicative, i.e. if $\gcd (m, n)=1$, then $\phi(mn)=\phi(m)\phi(n)$. We have Euler's product formula
 $$\phi(n)=n\prod_{p|n} \left(1-\frac{1}{p}\right),$$
and the following discrete Fourier representation \cite{Sch}, 
$$ \phi(n)=\sum_{k=1}^n \gcd (k, n) e^{-\frac{2\pi i k}{n}}.  $$

In order to estimate the size of the set $S^{\phi}_{M, N}$, we need to consider the number of solutions of the equation $\phi(n)=\phi(m)$. Given $m$, Carmichael (\cite{Carmi-1}, \cite{Carmi-2}) conjectured that there is at least one other integer $n\neq m$ such that $\phi(n)=\phi(m)$, which is the so called Carmichael's totient function conjecture. For each natural number $m$, let $A(m)$ be the number of $n$ such that $\phi(n)=m$. 
An alternative way of stating Carmichael's conjecture is that $A(m)$ can never be 1. 

We will use the bounds of $A(m)$ to verify Condition $(\star\star)$ for $\phi$. Ford \cite{Ford} showed that, for any $k\geq 2$, there exist infinitely many $m$ such that $A(m)=k$. For the upper bound, Pomerance \cite{P} showed that
\begin{equation}\nonumber
A(m)\leq m \exp\left(-(1+o(1)) \log m \log\log \log m/\log\log m \right)=: U(m).
\end{equation}
Therefore, 
\begin{equation}\label{S-phi}
|S^{\phi}_{M, N}|\lesssim M U(\phi(2M))\lesssim M U(M)\lesssim M\cdot \frac{M}{(\log M)^C} \text{~for any } C>0,  
\end{equation}
where the implied constant is absolute. Note that if $\phi(n)=m$, then $n \lesssim m \log\log m$ (\cite{Hardy-Wright}, Thereom 328). Hence $S^{\phi}_{M, N}$ is not empty only when $M\lesssim N\log\log N$ and $N\lesssim M\log\log M$. Combining this fact with \eqref{S-phi}, we get 
\begin{equation}\label{S-phi-delta}\nonumber
|S_{M,N}^{\phi}|\lesssim \frac{MN}{(\log M\log N)^{1+\delta}} \text{ for any } \delta>0,
\end{equation}
as desired. 

\subsection{Gaps between primes} Let $\pi(x)$ be the number of primes no more than $x$. The famous Prime Number Theorem (PNT) states that, as $x\rightarrow\infty$, 
\begin{equation}\label{PNT}
\pi(x)\sim \frac{x}{\log x}.
\end{equation}
Here $f(x)\sim g(x)$ means that $f(x)/g(x)$ goes to 1 as $x$ goes to $\nf$. PNT was first proved independently in 1896 by Hadamard and de la Vall\'{e}e Poussin. They both used the properties of the Riemann zeta-function $\zeta(s)=\sum_{n\geq 1} n^{-s}$ introduced by Riemann in his celebrated memoir.

Since $\pi(n)=\pi(m)$ can only occur when $n$ and $m$ are between two consecutive prime numbers, we need information about the gaps between primes, which have been extensively studied and many conjectures still remain open. PNT implies that the average gap between a prime $p$ and the next prime is about $\log p$. By \eqref{PNT}, one can derive that, for any $\epsilon>0$, there exists a prime in the interval $(p, p+\epsilon p]$ for sufficiently large prime $p$. However, this is not enough for our application, as we need results for primes in shorter intervals. 

Let $I(\theta, x)$ be the interval $[x, x+x^{\theta}]$. Hoheisel \cite{Hoh} showed that  $I(\theta, x)$ contains primes
for any $\theta>\frac{32999}{33000}$ as $x\rightarrow\infty$. 
Later, several authors made contributions to get smaller values of $\theta$ for which $I(\theta, x)$ contains primes for sufficiently large $x$. Iwaniec and Jutila (\cite{Iwa-Juti}, $\theta>\frac{5}{9}$) introduced to this problem a sieve method, which was later refined by Heath-Brown and Iwaniec (\cite{Heath-Iwa}, $\theta>\frac{11}{20}$).   The best result to date is  due to Baker, Harman, and Pintz \cite{Baker-Har-Pint}, who showed that we can take $\theta\geq 0.525$. Riemann Hypothesis implies that $I(\theta, x)$ contains primes for any $\theta>\frac{1}{2}$ as $x\rightarrow\infty$ \cite{Cramer}. 

Using the result of Baker-Harman-Pintz, the size of the set $S^{\pi}_{M, N}$
$$|S^{\pi}_{M, N}|\lesssim M^{\theta_0} M, $$
where $\theta_0=0.525$. By the above results about gaps between primes, we have $M\simeq N$ if $S^{\pi}_{M, N}\neq \emptyset$. Since $\theta_0<1$, we deduce that
\begin{equation}\nonumber
|S_{M,N}^{\pi}|\lesssim \frac{MN}{(\log M\log N)^{1+\delta}} \text{ for any } \delta>0,
\end{equation}
and hence we verify Condition $(\star\star)$ for $\pi$.

It is worth mentioning here some recent breakthroughs concerning gaps between primes. Let $p_n$ be the $n$-th prime. Zhang \cite{Zhang} and Maynard \cite{May} showed that there exists some absolute constant $C$ such that $p_{n+1}-p_n<C$ happens infinitely often. For large gaps, Ford, Green, Konyagin, Maynard, and Tao \cite{Ford-GKMT} proved that there are infinitely many $n$'s such that
$$
p_{n+1}-p_{n} \gtrsim \frac{\log n\log\log n\log\log\log\log n}{\log\log\log n}.
$$
It is possible that the boundedness of $B_{P,Q}^{\text{dis}}$ could provide a new approach to study gaps between primes. We shall not pursue this interesting idea here.

\section{Proof of Theorem \ref{new thm}} \label{section: proof}
\setcounter{equation}0

By symmetry, we only consider the case that $P$ satisfies Condition $(\star\star)$ with parameter $\delta$ and $\delta'$ (and $Q$ is arbitrary). For notational convenience, we will simply write $T$ for $B_{P,Q}^{\text{dis}}$. For any $\lambda>0$ and $f,g\in D(\ZZ)$, define the level set 
$$
E_{\lambda}:=\{n\in\ZZ: |T(f,g)(n)|>\lambda \}.
$$ 
Fix $\epsilon\in [\frac{1}{2\delta+1},1]$. Our goal is thus to prove the following the level set estimate
\begin{equation} \label{eq: goal}
|E_{\lambda}|\lesssim_{\delta,\delta'} \frac{1}{\lambda^{1+\epsilon}}, \text{ whenever } \|f\|_{l^2}=\|g\|_{l^2}=1.
\end{equation}
The (strong) $l^{1+\epsilon}$-bound of $T(f,g)$ will follow immediately from interpolation.

We will only consider the case $\lambda < 1$. The other case can be proved similarly (In fact, the case $\lambda\ge 1$ is simpler: just set $M=0$ in the proof below).

Define the Fourier transform for any $f\in D(\ZZ)$ by 
$$
\hat{f}(\xi):=\sum_{m\in\ZZ}f(m)e^{-2\pi i\xi m}.
$$ Then our operator can be rewritten as 
\[
\begin{split}
T(f,g)(n)&=\sum_{m\in\ZZ\setminus\{0\}} f(n-P(m))g(n-Q(m))\frac{1}{m}\\
&=\int_{\bT^2}\hat{f}(\xi)\hat{g}(\eta)e^{2\pi i(\xi+\eta)n}\sigma(\xi,\eta)\,d\xi d\eta,
\end{split}
\]
where $\bT$ is the unit circle and $\sigma$ is the periodic bilinear multiplier given by
\[
\sigma(\xi,\eta)=\sum_{m\in\ZZ\setminus\{0\}} \frac{1}{m}e^{-2\pi i(P(m)\xi+Q(m)\eta)}.
\]

Using a standard technique, we proceed to decompose the multiplier $\sigma$ dyadically. Choose an odd function $\rho\in\cs(\ZR)$ supported in the set $\{x:|x|\in(\frac{1}{2},2)\}$ with the property that
\[
\frac{1}{x}=\sum_{j=0}^{\nf}\frac{1}{2^j}\rho\left(\frac{x}{2^j}\right) \text{ for any } x\in\ZR \text{ with } |x|\ge 1.
\]
Let
\[
\sigma_j(\xi,\eta):=\frac{1}{2^j}\sum_{m\in\ZZ}\rho\left(\frac{m}{2^j}\right)e^{-2\pi i(P(m)\xi+Q(m)\eta)},
\]
and consequently $$\sigma(\xi,\eta)=\sum_{j=0}^{\nf} \sigma_j(\xi, \eta).$$

Correspondingly $T$ can be written as the sum $\sum_{j=0}^{\nf}T_j$, where
\[
\begin{split}
T_j(f,g)(n)&=\int_{\bT}\int_{\bT}\hat{f}(\xi)\hat{g}(\eta)e^{2\pi i(\xi+\eta)n}\sigma_j(\xi,\eta)\,d\xi d\eta \\
&=\frac{1}{2^j}\sum_{m\in\ZZ}\rho\left(\frac{m}{2^j}\right)f(n-P(m))g(n-Q(m)).
\end{split}
\]

Let $M$ be a non-negative integer to be determined later. Decompose $T$ into two parts: $\sum_{j=0}^{M-1}T_j$ and $\sum_{j=M}^{\nf}T_j$. It remains to control the level sets
\begin{equation}
E_{\lambda}^{(1)}:=\left\{n\in\ZZ: \left|\sum_{j=0}^{M-1}T_j(f,g)(n)\right|>\lambda \right\}
\end{equation}
and 
\begin{equation}
E_{\lambda}^{(2)}:=\left\{n\in\ZZ: \left|\sum_{j=M}^{\nf}T_j(f,g)(n)\right|>\lambda \right\}
\end{equation}

$E_{\lambda}^{(1)}$ can be estimated by the following simple lemma, whose proof is based on H\"older inequality and is omitted.  
\begin{lemma} \label{lemma}
There is an absolute positive constant $C$ such that for any $j\ge 0$, 
$$\|T_j(f,g)\|_{l^1}\le C\|f\|_{l^2}\|g\|_{l^2}.$$
\end{lemma}

By Lemma \ref{lemma} and triangle inequality, we see that
\begin{equation} \label{E1 estimate}
|E_{\lambda}^{(1)}|\le \frac{M}{\lambda}. 
\end{equation}
Note that \eqref{E1 estimate} is useful, i.e. better than the upper bound $\frac{1}{\lambda^{1+\epsilon }}$, only when $\lambda< 1$. We do not need this estimate in the case $\lambda\ge 1$.

To control $|E_{\lambda}^{(2)}|$, we employ a $TT^*$ method. Define an auxiliary function 
$$
h(n)=\frac{\overline{II(f,g)(n)}}{|II(f,g)(n)|}\hichi_{E_{\lambda}^{(2)}}(n),
$$ 
where
$$
II(f,g)(n):=\sum_{j=M}^{\nf}T_j(f,g)(n).
$$
It is easy to verify that
\begin{equation} \label{eq: TT* 1}
\lambda^2|E_{\lambda}^{(2)}|^2\le\left(\sum_{n\in\ZZ} II(f,g)(n)h(n)\right)^2.
\end{equation}
By Fubini's theorem and the definition of Fourier transform we obtain
\[
\begin{split}
\sum_{n\in\ZZ}II(f,g)(n)h(n)&=\int_{\bT^2}\hat{f}(\xi)\hat{g}(\eta)\sum_{j=M}^{\nf}\sigma_j(\xi,\eta)\hat{h}(-(\xi+\eta))\,d\xi d\eta \\
&=\int_{\bT^2}\hat{f}(\xi)\hat{g}(\eta-\xi)\sum_{j=M}^{\nf}\sigma_j(\xi,\eta-\xi)\hat{h}(-\eta)\,d\eta d\xi 
\end{split}
\]
Invoking Cauchy-Schwarz inequality and Plancherel's Theorem we get
\begin{equation} \label{eq: TT* 2}
\begin{split}
&\left(\sum_{n\in\ZZ}II(f,g)(n)h(n)\right)^2 \\
&\le \left(\int_{\bT} |\hat{f}(\xi)| |E_{\lambda}^{(2)}|^{\frac{1}{2}} \left(\int_{\bT} |\hat{g}(\eta-\xi)|^2\left|\sum_{j=M}^{\nf}\sigma_j(\xi,\eta-\xi) \right|^2 d\eta \right)^{\frac{1}{2}}d\xi \right)^2 \\
&\le |E_{\lambda}^{(2)}|  \int_{\bT} \int_{\bT} |\hat{g}(\eta)|^2\left|\sum_{j=M}^{\nf}\sigma_j(\xi,\eta) \right|^2 d\eta d\xi  \\
&\le V|E_{\lambda}^{(2)}|
\end{split}
\end{equation}
where
\[
V:=\sup_{\eta\in \bT}\int_{\bT} \left|\sum_{j=M}^{\nf}\sigma_j(\xi,\eta)\right|^2d\xi.
\]
Using \eqref{eq: TT* 1} and \eqref{eq: TT* 2}, we see that 
\begin{equation} \label{E2 square decay}
|E_{\lambda}^{(2)}|\le\frac{V}{\lambda^2}.
\end{equation}
To control $V$, we recall 
\[
S_{M,N}^P=\left\{(m,n): P(m)=P(n), \frac{1}{2}N\le |n|\le 2N, \frac{1}{2}M\le |m|\le 2M\right\}, 
\]
and note that for any $\eta\in\bT$,
\begin{equation} \label{eq: cancellation}
\begin{split}
&\int_{\bT}\left|\sum_{j=M}^{\nf}\sigma_j(\xi,\eta)\right|^2\,d\xi\\
&\le \sum_{j_1,j_2=M}^{\nf}\frac{1}{2^{j_1}}\frac{1}{2^{j_2}}\sum_{m_1,m_2\in\ZZ}\left|\rho\left(\frac{m_1}{2^{j_1}}\right)\rho\left(\frac{m_2}{2^{j_2}}\right)\right|\hichi_{S_{2^{j_1}, 2^{j_2}}^{P}}(m_1,m_2).
\end{split}
\end{equation}

By  the support of $\rho$ and Condition $(\star\star)$, \eqref{eq: cancellation} implies that
\begin{equation} \label{estimate for V}
V\lesssim \sum_{j_1,j_2=M}^{\nf}\frac{1}{2^{j_1}}\frac{1}{2^{j_2}}|S_{2^{j_1},2^{j_2}}^P|\lesssim_{\delta'} \sum_{j_1,j_2=M}^{\nf} \frac{1}{(j_1j_2)^{1+\delta}}\lesssim_{\delta, \delta'} \frac{1}{M^{2\delta}}.
\end{equation}
Combing \eqref{E2 square decay} and \eqref{estimate for V}, we get the estimate for $|E_{\lambda}^{(2)}|$
\begin{equation} \label{E2 estimate}
|E_{\lambda}^{(2)}|\lesssim_{\delta,\delta'} \frac{1}{\lambda^2M^{2\delta}}.
\end{equation}
Apply \eqref{E1 estimate} and \eqref{E2 estimate}, and one gets
\begin{equation}
|E_{\lambda}|\le |E_{\lambda/2}^{(1)}|+|E_{\lambda/2}^{(2)}|\lesssim_{\delta,\delta'} \frac{M}{\lambda}+\frac{1}{\lambda^2M^{2\delta}}.
\end{equation}
 
Optimize the above upper bound by choosing $M$ to be an integer near $(\frac{1}{\lambda})^{\frac{1}{1+2\delta}}$, and we obtain 
\[
|E_{\lambda}|\lesssim_{\delta,\delta'} \frac{M}{\lambda}\lesssim \left(\frac{1}{\lambda}\right)^{\frac{1}{1+2\delta}}\frac{1}{\lambda}\lesssim\frac{1}{\lambda^{1+\epsilon}},
\]
as $\epsilon\ge \frac{1}{1+2\delta}$ and $\lambda<1$.
This is our desired estimate \eqref{eq: goal}, and the proof of Theorem \ref{new thm} is thus complete.

\bigskip
\noindent \textbf{Acknowledgement.} The first author thanks Xiumin Du, Xiaochun Li and Ping Xi for many helpful discussions on related topics. The second author was partially supported by NSF grant DMS-1501982. Both authors thank the anonymous referee for valuable suggestions.

\end{document}